\documentclass[11pt]{article}

\usepackage[dvips]{graphics}
\usepackage{amsfonts}
\usepackage{amssymb}
\usepackage{amsthm}

\theoremstyle{plain}
\newtheorem{theorem}{Theorem}[section]
\newtheorem{lemma}[theorem]{Lemma}
\newtheorem{proposition}[theorem]{Proposition}
\newtheorem{corollary}[theorem]{Corollary}
\theoremstyle{definition}
\newtheorem{definition}[theorem]{Definition}

\newtheorem{remark}[theorem]{Remark}

\newtheorem{example}[theorem]{Example}

\newcommand{\xs}{x_1,\ldots,x_n}                %x_1,...,x_n
\newcommand{\vs}{v_1,\ldots,v_n}                %v_1,...,v_n
\newcommand{\Ms}{M_1,\ldots,M_q}                %M_1,...,M_q
\newcommand{\Fs}{F_1,\ldots,F_q}                %F_1,...,F_q

\newcommand{\depth}{{\rm{depth}}\ }              %depth
\newcommand{\mx}{{\rm{max}}\ }                  %max
\newcommand{\height}{{\rm{height}}\ }              %height
\newcommand{\dimn}{ {\rm{dim}} \ }
\renewcommand{\AA}{{\mathbb{A}}}                %Koszul complex A
\newcommand{\KK}{{\mathbb{K}}}                  %Koszul complex K 
\newcommand{\LL}{{\mathbb{L}}}                  %Koszul complex L
\newcommand{\F}{{\mathcal{F}}}                  %Facet ideal of a complex
\newcommand{\N}{{\mathcal{N}}}                  %Non-facet ideal of a complex
\newcommand{\U}{{\mathcal{U}}}                  %Universal facet of a leaf
\newcommand{\HH}{{\mathop{\rm{H}}}}             %Homology module
\newcommand{\df}{\delta_{\mathcal{F}}}          %facet  complex of an ideal
\newcommand{\dn}{\delta_{\mathcal{N}}}          %non-facet  complex of an ideal
\newcommand{\ndiv}{\not |}              
\newcommand{\st}{\ | \ }                        % such that: | 
\newcommand{\tuple}[1]{\langle #1 \rangle}      % < ... >
\title{\sc The Facet Ideal of a Simplicial Complex}
\author{Sara Faridi\thanks{Mathematics Department,
George Washington University,
Washington DC, 20052.
email: \emph{faridi@ gwu.edu.} %\newline
2000 Mathematics Subject classification: 13}}
\date{December 31, 2001} 
\begin{document}

\maketitle

\begin{abstract} To a simplicial complex, we associate a square-free
 monomial ideal in the polynomial ring generated by its vertex set
 over a field. We study algebraic properties of this ideal via
 combinatorial properties of the simplicial complex. By
 generalizing the notion of a tree from graphs to simplicial
 complexes, we show that ideals associated to trees satisfy sliding
 depth condition, and therefore have normal and Cohen-Macaulay Rees
 rings. We also discuss  connections with the theory of
 Stanley-Reisner rings.
\end{abstract}

          Given a graph on $n$ vertices, Villarreal ([Vi1]) defined the
         \emph{edge ideal} associated to that graph in a polynomial
         ring in $n$ variables (each variable representing a vertex of
         the graph) to be the ideal generated by monomials $xy$, where
         the corresponding vertices to $x$ and $y$ are connected by an
         edge. For example, the ideal $I=(xy, yu, yv, uv)$ corresponds
         to the following graph:

\[ \includegraphics{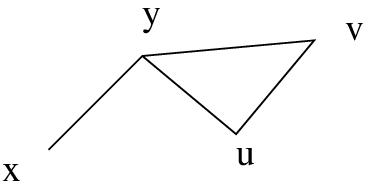} \]

        Later in [SVV], Simis, Vasconcelos and Villarreal used this
        construction along with properties of graphs to show that edge
        ideals of trees satisfy sliding depth condition. Among other
        things, this implies that the Rees ring of the edge ideal of a
        tree is normal and Cohen-Macaulay.

        Our goal here is generalize this construction to simplicial
        complexes. We define the notion of tree for simplicial
        complexes, and show that ideals corresponding to trees satisfy
        sliding depth (Theorem~\ref{main theorem}) and therefore have
        normal and Cohen-Macaulay Rees rings
        (Corollary~\ref{normal}). We also show that if the ideal of
        the tree is Cohen-Macaulay to begin with, it is strongly
        Cohen-Macaulay (Corollary~\ref{SCM}), meaning that all Koszul
        homology modules of generators of that ideal are
        Cohen-Macaulay.  Consequently we recover a rather large class
        of normal square-free monomial ideals with sliding depth
        condition.

        Traditionally, given a simplicial complex $\Delta$ one would
        associate to it the so-called \emph{Stanley Reisner} ideal,
        that is, the ideal generated by monomials corresponding to
        \emph{non-faces} of this complex (here again we are assigning
        to each vertex of the complex one variable of a polynomial
        ring generated by the vertices of the complex). For example,
        for the simplicial complex $\Delta$ below:
\[ \includegraphics{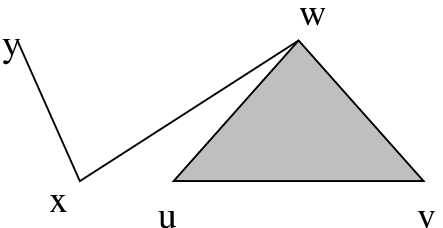} \]
the Stanley-Reisner ideal is $$\N(\Delta)=(xu, xv, yu, yv, yw).$$

        Our approach in this paper is to assign to the same simplicial
        complex $\Delta$ the ideal generated by its \emph{facets}:
        $$\F(\Delta)=(uvw, xw, xy).$$

\section{Basic Setup}

We first fix some notation and terminology.

\begin{definition}[simplicial complex, facet and more] 
A \emph{simplicial complex} $\Delta$ over a set of vertices $V=\{ \vs
\}$ is a collection of subsets of $V$, with the property that $\{ v_i
\} \in \Delta$ for all $i$, and if $F \in \Delta$ then all subsets of
$F$ are also in $\Delta$ (including the empty set). An element of
$\Delta$ is called a \emph{face} of $\Delta$, and the \emph{dimension}
of a face $F$ of $\Delta$ is defined as $|F| -1$, where $|F|$ is the
number of vertices of $F$.  The faces of dimensions 0 and 1 are called
\emph{vertices} and \emph{edges}, respectively, and $\dimn \emptyset
=-1$.

The maximal faces of $\Delta$ under inclusion are called \emph{facets}
of $\Delta$. The dimension of the simplicial complex $\Delta$ is the
maximal dimension of its facets; in other words $$\dimn \Delta =\max
\{ \dimn F \st F \in \Delta \}.$$ 

We denote the simplicial complex $\Delta$ with facets $\Fs$ by
$$\Delta = \tuple{\Fs}$$ and we call $\{ \Fs \}$ the \emph{facet set}
of $\Delta$. 

A simplicial complex with only one facet is called a \emph{simplex}.

\end{definition}

\begin{definition}[subcomplex]By a \emph{subcomplex} of a simplicial
 complex $\Delta$, in this paper, we mean a simplicial complex whose
facet set is a subset of the facet set of $\Delta$.
\end{definition}

\begin{definition}[facet ideal, non-face ideal] Let  $\Delta$
 be a simplicial complex over $n$ vertices labeled $\vs$. Let $k$ be a
field, $\xs$ be indeterminates, and $R$ be the polynomial ring
$k[\xs]$.

\noindent (a) We define $\F(\Delta)$ to be the ideal of $R$ generated
by square-free monomials $x_{i_1}\ldots x_{i_s}$, where
$\{v_{i_1},\ldots, v_{i_s}\}$ is a facet of $\Delta$. We call
$\F(\Delta)$ the \emph{facet ideal} of $\Delta$.

\noindent (b) We define $\N(\Delta)$ to be the ideal of $R$ generated
by square-free monomials $x_{i_1}\ldots x_{i_s}$, where $\{v_{i_1},\ldots,
v_{i_s}\}$ is not a face of $\Delta$. We call $\N(\Delta)$ the
\emph{non-face ideal} or the \emph{Stanley-Reisner ideal} of $\Delta$.
\end{definition}

        We refer the reader to [BH] for an extensive coverage of the
        theory of Stanley-Reisner ideals.

\begin{definition}[facet complex, non-face complex] Let
  $I=(\Ms)$ be an ideal in a polynomial ring $k[\xs]$, where $k$ is a
field and $\Ms$ are square-free monomials in $\xs$ that form a minimal
set of generators for $I$.

\noindent (a) We define $\df(I)$ to be the simplicial complex over a
set of vertices $\vs$ with facets $\Fs$, where for each $i$, $F_i=\{v_j
\st  x_j|M_i, \ 1 \leq j \leq n \}$. We call $\df(I)$ the \emph{facet
complex} of $I$.

\noindent (b) We define $\dn(I)$ to be the simplicial complex over a 
set of vertices $\vs$, where $\{v_{i_1},\ldots, v_{i_s}\}$ is a face
of $\dn(I)$ if and only if $x_{i_1}\ldots x_{i_s} \notin I$. We call
$\dn(I)$ the \emph{non-face complex} or the \emph{Stanley-Reisner
complex} of $I$.

\end{definition}

\begin{remark} It is worth observing that given a simplicial complex
 $\Delta$, if $$\F(\Delta)=(\Ms) \subseteq k[\xs],$$ then $\N(\Delta)$
 is generated by square free monomials that do not divide any of
 $\Ms$. \end{remark}

        We now generalize some notions from graph theory to simplicial
        complexes.

\begin{definition}[minimal vertex cover, vertex covering number] Let 
$\Delta$ be a simplicial complex with vertex set $V$ and facets
 $\Fs$. A \emph{vertex cover} for $\Delta$ is a subset $A$ of $V$, with
 the property that for every facet $F_i$ there is a vertex $v \in A$
 such that $v \in F_i$. A \emph{minimal vertex cover} of $\Delta$ is a
 subset $A$ of $V$ such that $A$ is a vertex cover, and no proper
 subset of $A$ is a vertex cover for $\Delta$. The smallest
 cardinality of a minimal vertex cover of $\Delta$ is called the
 \emph{vertex covering number} of $\Delta$.
\end{definition}

A simplicial complex may have several minimal vertex covers. 

\begin{example} Let $\Delta$ be the simplicial complex below.
\[ \includegraphics{stanleyexample.eps} \]

Here $\N(\Delta)=(xu, xv, yu, yv, yw)$, $\F(\Delta)=(uvw, xw, xy)$ and
$\dn \left (\F(\Delta) \right )$ is the following simplicial complex:
\[ \includegraphics{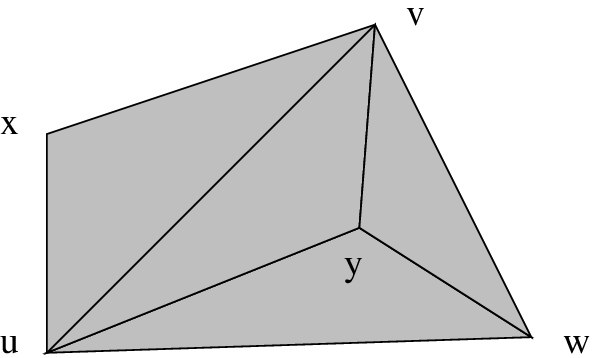} \]

Also, $\{x, w\}$, $\{y, w\}$, $\{x, v\}$, $\{x, u\}$ are all the
minimal vertex covers for $\Delta$.
\end{example}

\begin{proposition}\label{minprime} Let $\Delta$ be a simplicial
 complex over $n$ vertices. Consider the ideal $I=\F(\Delta)$ in the
 polynomial ring $k[\xs]$ over a field $k$. Then an ideal
 $p=(x_{i_1},\ldots,x_{i_s})$ of $R$ is a minimal prime of $I$ if and
 only if $\{x_{i_1},\ldots,x_{i_s}\}$ is a minimal vertex cover for
 $\Delta$.
\end{proposition}

        \begin{proof} Suppose that $I=(\Ms)$. The minimal primes of
         $I$ are generated by subsets of $\{\xs\}$ (see Proposition
         5.1.3 of [Vi2]). By definition $\{x_{i_1},\ldots,x_{i_s}\}$
         is a vertex cover for $\Delta$ if and only if for each
         generator $M_j$ of $I$, $x_{i_t} | M_j$ for some $1 \leq t
         \leq s$. It follows that $I \subseteq
         p=(x_{i_1},\ldots,x_{i_s})$ if and only if
         $\{x_{i_1},\ldots,x_{i_s}\}$ is a vertex cover for $\Delta$.
         The assertion now follows.  \end{proof}
         
\begin{definition}[unmixed] A simplicial complex $\Delta$ is 
\emph{unmixed} if all of its minimal vertex covers have the same
 cardinality.
\end{definition}

For instance, the simplicial complex in the previous example is unmixed.

\begin{definition}[pure] A simplicial complex $\Delta$ is
   \emph{pure} if all of its facets have the same
 dimension. Equivalently, this means that if $\F(\Delta)$ is generated
 by $\Ms$, all the $M_i$ are the product of the same number of
 variables.
\end{definition}

\begin{corollary}\label{pureiffunmixed} Let $I$ be a square-free 
monomial ideal in the polynomial ring $k[\xs]$. Then $\df(I)$ is
unmixed if and only if $\dn(I)$ is pure.\end{corollary}

        \begin{proof} By Proposition~\ref{minprime} $\df(I)$ is
        unmixed if and only if all minimal primes of $I$ have the same
        number of minimal generators, say that number is $s$. By the
        proof of Theorem~5.1.4 of [BH], $(x_{i_1},\ldots,x_{i_s})$ is
        a  minimal prime of $I$ if and only if $\{\vs \} \backslash
        \{ v_{i_1},\ldots,v_{i_s} \}$ is a facet of $\dn(I)$, which is
        equivalent to $\dn(I)$ being pure. \end{proof}

\begin{corollary}[A Cohen-Macaulay simplicial complex is unmixed] 
Suppose that $\Delta$ is a simplicial complex with vertex set
 $\xs$. If $k[\xs]/\F(\Delta)$ is Cohen-Macaulay, then $\Delta$
 is unmixed.\end{corollary}

        \begin{proof} If $k[\xs]/\F(\Delta)$ is Cohen-Macaulay, then
        $\dn(\F(\Delta))$ is pure ([BH] Corollary~5.1.5).
        Corollary~\ref{pureiffunmixed} then implies that
        $\Delta=\df(\F(\Delta))$ is unmixed. \end{proof}

\begin{remark}\label{dimension} It is worth observing that for a 
square-free monomial ideal $I$, there is a natural way to construct
$\dn(I)$ and $\df(I)$ from each other using the structure of the
minimal primes of $I$. To do this, consider the vertex set $V$
consisting of all variables that divide a monomial in the generating
set of $I$. The following correspondence holds:

        \begin{center} $F=$ facet of $\dn(I)$ $\longleftrightarrow $
        $V \setminus F=$ minimal vertex cover of $\df(I)$ \end{center}

Also, $$I =\bigcap p$$ where the intersection is taken over all
prime ideals $p$ of $k[V]$ that are generated by a minimal vertex
cover of $\df(I)$ (or equivalently, primes $p$ that are generated by
$V \setminus F$, where $F$ is a facet of $\dn(I)$; see [BH] Theorem 
5.1.4). 

        Regarding the dimension and codimension of $I$, note that by
        Theorem 5.1.4 of [BH] and the discussion above, setting
        $R=k[V]$ as above, we have
        $$\dimn R/I =\dimn \dn(I) +1 = |V| -{\rm \ vertex\ covering\
number\ of\ }\df(I)$$ and $$\height I\ ={\rm \ vertex\ covering\
number\ of\ } \df(I).$$

\end{remark}

We illustrate all this through an example.

\begin{example} For $I=(xy, xz)$,

\begin{center}
\begin{tabular}{cc}
 \underline{$\dn(I)$} & \underline{$\df(I)$} \\
&\\
$\includegraphics{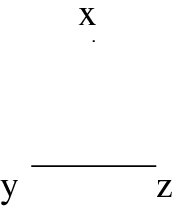}$ & $ \includegraphics{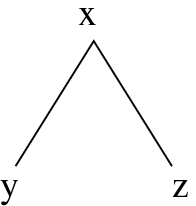}$ \\
&\\
\hspace{.3in}\underline{facets of $\dn(I)$}\hspace{.3in} &
\hspace{.3in}\underline{minimal vertex covers of
$\df(I)$}\hspace{.3in} \\ &\\ $\{ x \}$ & $\{ y, z \}$ \\ $\{ y, z
\}$& $\{ x \}$
\end{tabular}
\end{center}
 Note that $I = (x) \cap (y,z)$, and $$\dimn k[x,y,z]/(xy,xz)=2$$
 as  asserted in Remark~\ref{dimension} above.

\end{example}

%%%%%%%%%%%%%%%%%%%%%%%%%%%%%%%%%%%%%%%%%%%%%%%%%%%%%%%%%%%%%%%%%%%%%%%%%%%%%

\section{Simplicial complexes that are trees}

         We now explore a new notion of \emph{tree} on simplicial
        complexes.  This definition generalizes trees in graph theory,
        and turns out to behave well under localization and removal of
        a facet as described below.

        To motivate the definition, recall that a connected graph is a
         tree if it has no cycles. An equivalent definition states
         that a connected graph is a tree if every subgraph has a
         \emph{leaf}, where a leaf is a vertex that belongs to only
         one edge of the graph. This latter description is the one
         that we adapt, with a slight change in the definition of a
         leaf, to the class of simplicial complexes.

\begin{definition}[leaf, universal set]~\label{leaf} Suppose that $\Delta$
 is a simplicial complex. A facet $F$ of $\Delta$ is called a
\emph{leaf} if either $F$ is the only facet of $\Delta$, or there
exists a facet $G$ in $\Delta$, $G \neq F$, such that $$ F \cap F'
\subseteq F \cap G$$ for every facet $F' \in \Delta$, $F' \neq F$.

We denote the set of facets $G$ in $\Delta$ with this property by
$\U_\Delta (F)$ and call it the \emph{universal set} of $F$ in
$\Delta$.

\end{definition}

Note that the facet $G$ in Definition~\ref{leaf} is not necessarily unique:

\begin{example} In the following simplicial complex $\Delta$, $F$ is a leaf
 and both $G_1$ and $G_2$ are in $\U_\Delta (F)$.

\[ \includegraphics{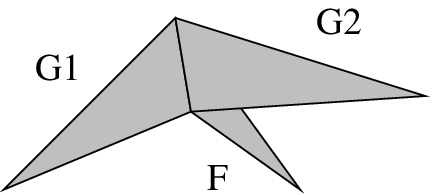} \]

\end{example}

\begin{remark}\label{vertex} In order to be
 able to quickly identify a leaf in a simplicial complex, it is
 important to notice that a leaf must have a vertex that does not
 belong to any other facet of that simplicial complex. This follows
 easily from Definition~\ref{leaf}: otherwise, a leaf $F$ will be
 contained in the members of its universal set, which contradicts the
 fact that a leaf is a facet. \end{remark}

\begin{definition}[tree]~\label{tree} Suppose that $\Delta$ 
is a connected simplicial complex.  We say that $\Delta$ is a
 \emph{tree} if every nonempty subcomplex of $\Delta$ (including
 $\Delta$ itself) has a leaf. 

Equivalently, $\Delta$ is a tree if every nonempty
 \emph{connected} subcomplex of $\Delta$ has a leaf.

Recall that by a subcomplex of $\Delta$ we mean a simplicial complex
whose facet set is a subset of the facet set of
$\Delta$. \end{definition}

\begin{definition}[forest]\label{forest} A simplicial complex
 $\Delta$ with the property that every connected component of $\Delta$
 is a tree is called a \emph{forest}. In other words, a forest is a
 simplicial complex with the property that every nonempty subcomplex
 has a leaf.
\end{definition}

\begin{example} The simplicial complex below is a tree.

 \[ \includegraphics{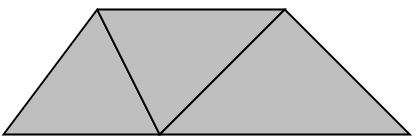} \]
\end{example}

\begin{example} The simplicial complex below is not a tree because every
 vertex is shared by at least two facets (see Remark~\ref{vertex}
 above).

\[ \includegraphics{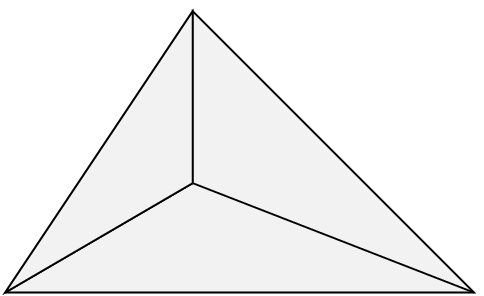} \]
\end{example}

\begin{example} The simplicial complex below with facets $F_1=\{a,b,c \}$,
 $F_2=\{ a,c,d \}$ and $F_3=\{ b,c,d,e \}$ is not a tree because the
only candidate for a leaf is the facet $F_3$, but neither one of $F_1
\cap F_3$ or $F_2 \cap F_3$ is contained in the other.

\[ \includegraphics{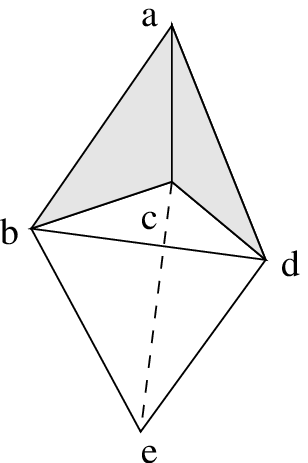} \]
\end{example}

        Notice that in the case that $\Delta$ is a graph,
        Definition~\ref{tree} agrees with the definition of a tree in
        graph theory, with the difference that now the term ``leaf''
        refers to an edge, rather than a vertex.

        A notion that will be crucial in the following discussion is
        ``removing a facet'' from a simplicial complex. We therefore
        record a definition for this construction.

\begin{definition}[facet removal]\label{removal} Suppose $\Delta$
 is a simplicial complex with facets $\Fs$ and $\F(\Delta)=(\Ms)$ its
 ideal in $R=k[\xs]$. The simplicial complex obtained by
 \emph{removing the facet} $F_i$ from $\Delta$ is the simplicial
 complex
 $$\tuple{F_1,\ldots,\hat{F}_{i},\ldots,F_q}.$$ \end{definition}

        An important property of a tree that we will use later is that
        the localization of a tree is a forest. Below, by abuse of
        notation, we use $\xs$ to denote both the vertices of a
        simplicial complex and the variables of the polynomial ring
        corresponding to that complex.

\begin{lemma}[Localization of a tree is a forest]\label{localization}
 Let $\Delta$ be a tree with vertices $\xs$, and let $I= \F (\Delta)$
 be the facet ideal of $\Delta$ in the polynomial ring $R=k[\xs]$
 where $k$ is a field. Then for any prime ideal $p$ of $R$, $\df
 (I_p)$ is a forest. \end{lemma}

        \begin{proof} The first step is to show that it is enough to
        prove this for prime ideals of $R$ generated by a subset of
        $\{ \xs \}$. To see this, assume that $p$ is a prime ideal of
        $R$ and that $p'$ is another prime of $R$ generated by all
        $x_i \in \{ \xs \}$ such that $x_i \in p$ (recall that the
        minimal primes of $I$ are generated by subsets of $\{ \xs
        \}$). So $p' \subseteq p$. If $I=(\Ms)$, then $$I_{p'} =
        ({M_1}',\ldots,{M_q}')$$ where for each $i$, ${M_i}'$ is the
        image of $M_i$ in $I_{p'}$. In other words, ${M_i}'$ is
        obtained by dividing $M_i$ by the product of all the $x_j$
        such that $x_j| M_i$ and $x_j \notin p'$. But $x_j \notin p'$
        implies that $x_j \notin p$, and so it follows that ${M_i}'
        \in I_p$. Therefore $I_{p'} \subseteq I_p$. On the other hand
        since $p' \subseteq p$, $I_p \subseteq I_{p'}$, which implies
        that $I_{p'} = I_p$ (the equality and inclusions of the ideals
        here mean equality and inclusion of their generating sets).

        We now prove the theorem for $p=(x_{i_1},\ldots, x_{i_r})$.
        The main point to notice here is that the simplicial complex
        corresponding to $I_p$ is obtained by removing all vertices
        except for $x_{i_1},\ldots, x_{i_r}$ from $\Delta$.  So the
        proof reduces to showing that if we remove a vertex from
        $\Delta$, the resulting simplicial complex is a forest.
        
        Let $\Delta= \tuple{  \Fs}$, and suppose that we remove the vertex
        $x_1$. Setting $${F_i} ' =F_i \setminus \{ x_1 \}$$ we would
        like to show that $$\Delta ' = \tuple{ F_1 ',\ldots, F_q '}$$ is a
        forest. Notice that some of the $F_i'$ appearing in $\Delta '$
        may be redundant. We need to show that every nonempty
        subcomplex of $\Delta '$ has a leaf.

         Let $\Delta_1 '=\tuple{ F_{i_1} ',\ldots, F_{i_s} '}$ be a
         subcomplex of $\Delta '$ where $F_{i_1} ',\ldots, F_{i_s} '$
         are distinct facets. Consider the corresponding subcomplex
         $\Delta_1=\tuple{ F_{i_1} ,\ldots, F_{i_s} }$ of $\Delta$,
         which has a leaf, say $F_{i_1}$. So there exists $G \in
         \Delta_1$ such that $$ F_{i_1} \cap F \subseteq F_{i_1} \cap
         G$$ for every facet $F \in \tuple{F_{i_2} ,\ldots, F_{i_s}
         }$. Now since each of the $ F_{i_j}'$ is a nonempty facet of
         $\Delta_1 '$ the same statement holds in $\Delta_1 '$, that
         is,
        $$ F_{i_1}' \cap F' \subseteq F_{i_1}' \cap G'$$ for every
        facet $F' \in \Delta_1' \setminus \{ F_{i_1}' \}$.
        This implies that $ F_{i_1}'$ is a leaf for $\Delta_1 '$.

        \end{proof}

%%%%%%%%%%%%%%%%%%%%%%%%%%%%%%%%%%%%%%%%%%%%%%%%%%%%%%%%%%%%%%%%%%%%%%%
      
\section{Cohen-Macaulay properties of simplicial complexes}

        The main result of this section states that the facet ideal of
        a tree (see Definition~\ref{tree}) has sliding depth in the
        polynomial ring generated by its set of vertices over a field
        $k$. A special case of this result for the case that $\Delta$
        is a graph was proved in [SVV] (Theorem~1.3). We first define
        the notion of sliding depth. Then we carry out a similar
        argument to the one in [SVV] to show that the facet ideal of a
        tree has sliding depth.

\begin{definition}[sliding depth] Let $I$ be an ideal in a ring
 $R$ of dimension $n$.  Let $\KK$ be the Koszul complex on the set
 ${\bf a}= \{a_1,\ldots,a_q\}$ of generators of $I$, and denote by
 $\HH_i(\KK)$ the homology modules of $\KK$. $I$ satisfies
 \emph{sliding depth} if
 $$\depth \HH_i(\KK) \geq n-q+i,  \hspace{.2 in} 
{\rm{\ for\  all \ }} i.$$
\end{definition}

\begin{definition}[Cohen-Macaulay simplicial complex] A simplicial complex 
$\Delta$ over a set of vertices labeled $\xs$ is Cohen-Macaulay, if
the quotient ring $k[\xs]/\F(\Delta)$ is
Cohen-Macaulay. \end{definition}

        Below, we assume that the polynomial ring $k[\xs]$ is
        localized at the maximal ideal $(\xs)$. By abuse of notation,
        we use $\xs$ to denote both the vertices of a simplicial
        complex and the variables of the polynomial ring corresponding
        to that complex.

\begin{theorem}[Main Theorem: trees have sliding depth]\label{main theorem} 
If a simplicial complex $\Delta$ on a set of vertices $\{ \xs
 \}$ is a tree and $k$ is a field, then $\F(\Delta)$ has sliding depth
 in the polynomial ring $R=k[\xs]$.
\end{theorem}

        \noindent{\bf Proof of Main Theorem.} Suppose that $\Delta$
        is a tree with facets $\Fs$ and vertex set $\xs$, and that 
        $\F(\Delta)= (\Ms)$, where each $M_i$ is a square-free
        monomial in the variables $\xs$. We want to show that
        $\F(\Delta)$ has sliding depth in the polynomial ring
        $R=k[\xs]$. We argue by induction on $q$.

        The case $q=1$ is the case where $\F(\Delta)=(M_1)$, and the
        Koszul complex looks like
        $$0 \longrightarrow R \stackrel{.M_1}{\longrightarrow}
        R\longrightarrow0 \ .$$ 
        which gives $\HH_0(M_1) = R/(M_1)$, and 
        $\depth R/(M_1) =n-1$, and so $\F(\Delta)$ has sliding depth.

        Suppose we know that the theorem holds for any tree with up to
        $q-1$ facets, $q \geq 2$. We want to show this for a
        simplicial complex with $q$ facets.

        Suppose without loss of generality that $F_q$ is a leaf and
        $$M_q=x_1 \ldots x_r.$$ By removing $F_q$ from $\Delta$ 
        (see Definition~\ref{removal}) we
        obtain a tree $$\Delta'=\tuple{F_1,\ldots, F_{q-1}}.$$

        We let $$Z' \subseteq \{\xs \}$$ be the vertex set for
         $\Delta'$, and $$Z''= \{ x_1, \ldots, x_r \} \setminus Z'.$$
         It follows that
          $$Z' \cap Z ''= \emptyset \ \
        {\rm{and}}\ \ \ Z' \cup Z'' =\{\xs\}.$$

         By the induction hypothesis, $\F(\Delta')$
        has sliding depth in $R'=k[Z']$.

        Let $y'$ and $y''$ be the product of the elements of $Z' \cap
        \{ x_1, \ldots, x_r \}$ and $Z''$, respectively, so that
        $$M_q=x_1 \ldots x_r= y'y''.$$  Let $$L'= \F(\Delta') + (y').$$ 
        
         Since $y'$ is the monomial describing  the
         intersection of the leaf $F_q$ with $\Delta'$, there is an
         $M_j$ in the generating set of $\F (\Delta')$ such that $y' |
         M_j$, and so adding $y'$ to $\F (\Delta ')$ does not increase
         the number of generators of $\F (\Delta ')$.

         Notice that $\df(L')$ is a forest with at most $q-1$
        facets. To show that it is a forest, we must show that every
        connected subcomplex of $\df(L')$ has a leaf. Let $F'$ denote
        the facet in $\df(L')$ corresponding to the monomial $y'$. The
        only subcomplexes that we need to worry about are those that
        contain $F'$, since any other one will be a connected
        subcomplex of $\Delta$, and so will have a leaf by definition.

        Suppose that $$\Delta''=\tuple{F', F_{i_1}, \ldots, F_{i_s}}$$ is a
        connected subcomplex of $\df(L')$. Consider the connected
        subcomplex of $\Delta$ defined as
        $$\overline{\Delta}=\tuple{F_q, F_{i_1}, \ldots, F_{i_s}}.$$ Since
        $\Delta$ is a tree, $\overline{\Delta}$ has a leaf. There are two
        possibilities:

        \noindent \emph{Case 1.} For some $j$, $F_{i_j}$ is a leaf of
         $\overline{\Delta}$, and $G \in
         \U_{\overline{\Delta}}(F_{i_j})$. Then $$ F_q \cap F_{i_j}
         \subseteq G \cap F_{i_j}.$$ But $F_q \cap F_{i_j} = F' \cap
         F_{i_j}$ (from the construction above), and so $ F_{i_j}$ is
         a leaf for $\Delta''$.

        \noindent \emph{Case 2.} $F_q$ is a leaf of
        $\overline{\Delta}$, and $G \in
        \U_{\overline{\Delta}}(F_q)$. Then for all $j$, $$F_{i_j} \cap
        F_q\subseteq G \cap F_q$$ which, as above, translates into
        $$F_{i_j} \cap F' \subseteq G \cap F'$$ which implies that
        $F'$ is a leaf for $\Delta''$.

        Therefore $\df{(L')}$ is a forest.  By the induction hypothesis,
        the facet ideal of every connected component of $\df(L')$
        satisfies sliding depth, and so $L'$ will satisfy sliding
        depth (this is because the Koszul homology of $L'$ can be
        written as a direct sum of the tensor products, over $k$, of
        the Koszul homologies of the facet ideal of each connected
        component of $\df(L')$; see the discussion on page 397 of [SVV]).
        
         If $R'=k[Z']$ and $R'' = k[Z'']$, then $$R=R' \otimes_k R''$$ and
         $$\F(\Delta)= \F(\Delta')+ (y)$$ where
         $$y=M_q=y'y''.$$
        
        If $\KK'$ denotes the Koszul complex of $\F(\Delta')$ over
        $R'$, $|Z'|=n'$ and $|Z''|=n''$, the sliding depth condition on
        $\F(\Delta')$ implies that for all $i$: $$\depth \HH_i(\KK')
        \geq n'-(q-1)+i = n'-q +i+ 1. $$

        If $\KK$ is the Koszul complex of $\F(\Delta')$ in
        $R=R'[Z'']$, since $R$ is a faithfully flat extension of $R'$
        by $n''$ variables, it follows from [BH] Proposition~1.6.7 and
        the inequality above that for every $u$,

$$\begin{array}{rl}
\depth \HH_i(\KK)&=\depth \HH_i(\KK') + \depth k[Z'']\\
& \geq n'- q + (i+1) + n''\\
& = n-q + (i+1).
\end{array}$$
        
        Suppose that $\LL$ is the Koszul complex of
        $\F(\Delta)=\F(\Delta')+ (y)$. As in Proposition 3.7 of [SVV],
        the basic relationship between the homologies of $\LL$ and
        $\KK$ is: \begin{eqnarray}\label{ses10} 0 \longrightarrow
        \overline{\HH_i(\KK)} \longrightarrow \HH_i(\LL)
        \longrightarrow {_y}\HH_{i-1}(\KK) \longrightarrow 0,
        \end{eqnarray} where for any module $E$, $\overline{E}=E/Ey$
        and $_yE=\{e\in E \st y.e=0\}$.

        We will need the following lemma:

\begin{lemma}[Depth Lemma]\label{depth lemma} If  $0\longrightarrow N 
\longrightarrow M \longrightarrow L \longrightarrow 0$ is a short
exact sequence of modules over a local ring $R$, then

(a) If $ \depth M < \depth L$, then $ \depth N= \depth M$.

(b) If $\depth M =\depth L $, then $\depth N \geq \depth M$.

(c) If $\depth M > \depth L $, then $\depth N= \depth L +1$.
\end{lemma}

        \begin{proof} This is Lemma 1.3.9 in [Vi2]. See Corollary 18.6
        in [E] for a proof. \end{proof}

\begin{lemma}\label{the lemma} With $\KK'$ defined (as above) as the
 Koszul complex of $\F(\Delta')$ over the polynomial ring $R'=k[Z']$,
 we have:
$$ \depth \overline{\HH_i(\KK')} \geq n'-q+i+1 $$
$$ \depth _{y'}\HH_{i}(\KK') \geq n'-q+i+1$$
where $\overline{\HH_i(\KK')}= \HH_i(\KK')/y' \HH_i(\KK')$, and  
$_{y'}\HH_i(\KK')=\{ x \in \HH_i(\KK') \st y'.x=0 \}$.
 \end{lemma}

        \begin{proof} Suppose that $$L'=\F(\Delta')+(y').$$ Notice
        that $L'$ is the facet ideal of a simplicial complex with $n'$
        vertices and $q'$ facets, where $q'\leq q-1$, since $y'$ is
        part of at least one facet of $\Delta'$. By induction $L'$ has
        sliding depth, so if $\AA'$ is the Koszul complex of $L'$,
        $$\depth \HH_i(\AA') \geq n'-q'+i.$$

        We also have an exact sequence: 
        $$ 0 \longrightarrow {_{y'}}\HH_{i}(\KK') \longrightarrow
        \HH_{i}(\KK')\stackrel{y'}{\longrightarrow} \HH_{i}(\KK')
        \longrightarrow \overline{\HH_i(\KK')} \longrightarrow 0$$
        which breaks into two exact sequences:
        \begin{eqnarray}\label{ses1} 0 \longrightarrow
        {y'}\HH_{i}(\KK') \longrightarrow \HH_{i}(\KK')
        \longrightarrow \overline{\HH_i(\KK')} \longrightarrow 0
        \end{eqnarray}
        \begin{eqnarray}\label{ses3} 0 \longrightarrow
        {_{y'}}\HH_{i}(\KK') \longrightarrow \HH_{i}(\KK')
        \longrightarrow {y'} \HH_i(\KK') \longrightarrow 0 \
        .\end{eqnarray}
        
        Similar to  Sequence~(\ref{ses10}) above, we also have the
        following short exact sequence: \begin{eqnarray}\label{ses2} 0
        \longrightarrow \overline{\HH_i(\KK')} \longrightarrow
        \HH_i(\AA') \longrightarrow {_{y'}}\HH_{i-1}(\KK')
        \longrightarrow 0 \ . \end{eqnarray}

        We know that $\KK'$ and $\AA'$ have sliding depth. We prove
        the lemma  by induction on $i$.

        For $i=0$, $\HH_{0}(\KK')=R'/{\F(\Delta')}$ has depth $\geq
        n'-q+1$. Since the facet corresponding to the monomial $y'$
         is contained in a facet of $\Delta'$,
        $$\overline{\HH_{0}(\KK')}=R'/(\F(\Delta'), y')$$ is the zeroth
        homology of a forest with at most $q-1$ facets, 
        and so in $R'$
        $$\depth \overline{\HH_{0}(\KK')} \geq n'-(q-1).$$

        Plugging this into Sequence~(\ref{ses1}), along with
        Lemma~\ref{depth lemma} implies that $$\depth
        {y'}\HH_{0}(\KK') \geq n' - (q -1)$$ and then
        Sequence~(\ref{ses3}) and Lemma~\ref{depth lemma} imply that:
        $$\depth _{y'}\HH_{0}(\KK') \geq n' - (q-1).$$

        Now suppose that the statement holds for all $j \leq i-1$. We
        verify it for $j=i$.

        Since $$\depth {_{y'}}\HH_{i-1}(\KK') \geq n'- (q-1)+i-1$$ and
        $$ \depth \HH_{i}(\AA') \geq n'- (q-1)+i,$$
        Sequence~(\ref{ses2}) and Lemma~\ref{depth lemma} imply that
        $$\depth \overline{\HH_i(\KK')} \geq n'- (q-1)+i.$$

        We plug this into Sequence~(\ref{ses1}) and consider the fact
        that $ \depth \HH_{i}(\KK') \geq n'- (q-1)+i$ along with
        Lemma~\ref{depth lemma} and conclude that $$ \depth
        {y'}\HH_{i}(\KK') \geq n'- (q-1)+i$$ which when plugged in
        Sequence~(\ref{ses3}) implies that $$ \depth
        {_{y'}}\HH_i(\KK') \geq n'- (q-1)+i.$$
        \end{proof}

\begin{proposition}\label{final step} If $\KK$ is the Koszul complex of $\F(\Delta')$ in $R$ and  $y=M_q=y'y''$, as above, then 
$$ \depth \overline{\HH_i(\KK)} \geq n-q+i$$ 
$$ \depth _{y}\HH_{i}(\KK) \geq n-q+i+1.$$
\end{proposition}

        \begin{proof} Notice that
        $$\HH_i(\KK)=\HH_i(\KK')\otimes _k k[Z'']$$ where $\KK'$
         and $Z''$ are defined above,  and
        multiplication by $y=y'y''$ in $\HH_i(\KK)$ can 
        be factored as
        $$ \HH_i(\KK')\otimes _k k[Z''] 
        \stackrel{y'\otimes_k
        1}{\longrightarrow} \HH_i(\KK')\otimes _k k[Z'']
        \stackrel{1\otimes_k y''}{\longrightarrow}\HH_i(\KK')\otimes _k
        k[Z''].$$
        
        From here, setting $E= \HH_i(\KK)$, $E'=\HH_i(\KK')$,
        $E''=k[Z'']$ and $y=y'y''$ we get an exact
        sequence:
 \begin{eqnarray}\label{ses4} 0\rightarrow {_{y'}}E'\otimes_k E''
\rightarrow {_{y}}E \rightarrow E' \otimes_k {_{y''}}E''\rightarrow
E'/y'E' \otimes_k E'' \stackrel{f}{\rightarrow} E/yE \rightarrow E'
\otimes_k E''/y''E'' \rightarrow 0 \end{eqnarray}

        The map $f$ is multiplication by $1\otimes_k y''$ and so the
        image of $f$ is isomorphic to
        $$E'/y'E' \otimes_k y''E''$$ which by Lemma~\ref{the lemma}
        has depth at least $n-q+(i+1)$, and 
        $$\depth E' \otimes_k
        E''/y''E'' \geq n'- (q-1)+
        i + n''-1 = n-q+i,$$
        and therefore by Depth Lemma:
        $$\depth E/yE \geq n-q+i.$$

        Also, since $ E''$ is a polynomial ring,
        ${_{y''}}E''=0$ and so
        $$\depth {_{y}}E= \depth {_{y'}}E'\otimes_k E'' =
        n-q+(i+1)$$ as desired.

         \end{proof}

\noindent{ \bf Conclusion of Proof of Main Theorem.} We now apply the
result of Proposition~\ref{final step} and Depth Lemma to the short 
exact sequence~(\ref{ses10}) from above:
        $$0 \longrightarrow \overline{\HH_i(\KK)} \longrightarrow
\HH_i(\LL) \longrightarrow {_y\HH_{i-1}(\KK)} \longrightarrow 0$$

and conclude that for all $i$,
$$\depth \HH_i(\LL) \geq n-q+i.$$
\hfill $\square$

\begin{proposition}\label{inequality} Suppose that $\Delta$ is a 
Cohen-Macaulay tree and $I= \F(\Delta) \subseteq k[\xs]$. Let $p$ be
 a prime ideal of $R$. Then, if $\mu (J)$ denotes the minimal number
 of generators of the ideal $J$, $$\mu (I_p) \leq \mx \{ \height I,
 \height p -1 \}.$$
\end{proposition}

\begin{proof}  If $p$ is a minimal prime of $I$, then
by Proposition~\ref{minprime} $p$ is generated by a minimal vertex
cover of $\Delta$, say $$p= (x_{i_1}, \ldots, x_{i_r}).$$ If
$I=(\Ms)$, then since $\{x_{i_1}, \ldots, x_{i_r}\}$ is a minimal
vertex cover, for each $x_{i_e} \in \{x_{i_1}, \ldots, x_{i_r}\}$
there is at least one $M_j \in \{\Ms \}$ such that $x_{i_e} | M_j$ and
$x_{f} \ndiv M_j$ for $x_f \in \{ x_{i_1}, \ldots, x_{i_r} \}
\setminus \{ x_{i_e} \}$. This implies that $p \subseteq I_p$ and
since $I_p \subseteq p$, it follows that $I_p =p$.

        So in particular, $$\mu (I_p) = \mu (p) = \height I$$ and
        therefore the theorem holds.

        Now we prove the general version of the theorem by induction
        on the number of vertices of $\Delta$. The case $n=1$ follows
        from the argument above. Suppose that the inequality holds for
        any simplicial complex with less than $n$ vertices. 

        Let $\Delta$ be a simplicial complex with vertices $\xs$, and
        let $p$ be a prime ideal of $k[\xs]$. We know by
        Lemma~\ref{localization} that $\df(I_p)$ is a forest. By the
        first paragraph of the proof of the same lemma, we can assume
        that $p$ is generated by a subset of $\{ \xs \}$.

        Now suppose that $p$ is not a minimal vertex cover. So, with
        $s > 0$, we have $$I_p = I_1 + \ldots + I_s + (x_{j_1},
        \ldots, x_{j_t})$$ where $\df (I_i)$ for $ i=1, \ldots, s$ and
        $\df(x_{j_e})$ for $e=1,\ldots,t$, are the connected
        components of the forest $\df(I_p)$. 

        Notice that each $\df(I_i)$ is also a Cohen-Macaulay tree, in
        the sense that if $R_i$ denotes the polynomial ring over $k$
        generated by the variables appearing in $I_i$, then $R_i/I_i$
        is a Cohen-Macaulay ring. This follows, for example, from [BH]
        Theorem~2.1.7, since $R_p/I_p$, which is Cohen-Macaulay, is
        isomorphic to the tensor product, over $k$, of all the
        $R_i/I_i$.

        Let $$p= p_1 +\ldots + p_s + (x_{j_1}, \ldots, x_{j_t})$$
        where each $p_i$ is the ideal of all vertices of $\df(I_i)$
        above.
        
        Since $p_i$ is the ideal generated by all vertices of
        $\df(I_i)$, it \emph{properly} contains the ideal generated by
        a minimal vertex cover of $\df(I_i)$. Therefore, by the first
        paragraph of this proof $$\height I_i < \height p_i$$ which
        implies that $$\height I_i \leq \height p_i -1.$$

        By the induction hypothesis for each $i$, since $I_i
        =(I_i)_{p_i}$,
        $$\mu(I_i) \leq \mx \{ \height I_i, \height p_i -1 \} \leq
        \height p_i -1.$$

        We now have (recall that $s > 0$)
        $$\begin{array}{rl} 
        \mu(I_p)&= \mu(I_1) + \ldots + \mu(I_s) + t \\
                &\leq \height p_1 -1 + \ldots + \height p_s -1 + t \\
                &= \height p -s \\
                & \leq \height p - 1
        \end{array}$$   

        The assertion then follows.
          \end{proof}

\begin{definition}[Strongly Cohen-Macaulay] An ideal $I$ 
of a ring $R$ is \emph{strongly Cohen-Macaulay} if all Koszul homology
modules of $I$ are Cohen-Macaulay. \end{definition}

\begin{corollary}[The facet ideal of a C-M tree is strongly C-M]\label{SCM}
 Suppose that $\Delta$ is a Cohen-Macaulay tree. Then $\F (\Delta)$ is
 a strongly Cohen-Macaulay ideal. \end{corollary}

        \begin{proof} This follows from Proposition~\ref{inequality}
        above, the fact that $\F (\Delta)$ has sliding depth and
        Theorem~1.4 of [HVV], or Theorem 3.3.17 of [V]. \end{proof}

\begin{definition}[{[V]}] An ideal $I$ of $R$ satisfies condition $\F _1$ if 
$\mu (I) \leq \dimn R_p$ for all prime ideals of $p$ of $R$ such that
$I \subset p$. \end{definition}

\begin{proposition}\label{F1} Suppose that $\Delta$ is a tree
 and $I= \F(\Delta) \subseteq k[\xs]$. Then $I$ satisfies $\F_1$.
\end{proposition}

        \begin{proof} As in Proposition~\ref{inequality}, if $p$ is a
         minimal prime of $I$ we have $$\mu (I_p) = \mu (p) = ht (p) =
         \dimn R_p. $$

        If $I$ represents a tree with $q$ facets, then we proceed by
        induction on $q$. The case $q=1$ results in $\mu (I_p) = 1
        \leq \dimn R_p$ for all primes $p$ containing $I$. Suppose
        that the proposition holds for all trees that have less than $q$
        facets, and let $I$ be the facet ideal of a tree with $q$
        facets, and $p$ be a prime ideal of $R$ that contains $I$. If
        $p$ is not minimal over $R$, similar to the proof of
        Proposition~\ref{inequality}, we have that $$I_p = I_1 +
        \ldots + I_s + (x_{j_1}, \ldots, x_{j_t})$$ where $\df (I_i)$
        for $ i=1, \ldots, s$ and $\df(x_{j_e})$ for $e=1,\ldots,t$,
        are the connected components of the forest $\df(I_p)$.

        Let $$p= p_1 +\ldots + p_s + (x_{j_1}, \ldots, x_{j_t})$$
        where each $p_i$ is the ideal of all vertices of $\df(I_i)$
        above, so for $i=1, \ldots, s$, $I_i=(I_i)_{p_i}$. By
        induction hypothesis, $ \mu(I_i) \leq \dimn R_{p_i}$ for all
        $i$, and so $$\mu (I_p) =\mu ( I_1) + \ldots + \mu (I_s) + t
        \leq \dimn R_{p_1} + \ldots + \dimn R_{p_s}+ t=\dimn R_p.$$

        \end{proof}

\begin{corollary}[The Rees ring of a tree is normal and C-M]\label{normal}
 Suppose that $\Delta$ is a tree over vertices $\xs$ with facet ideal
 $I=\F (\Delta)$ in the polynomial ring $k[\xs]$ where $k$ is a
 field. Then the Rees ring of $I$ is normal and
 Cohen-Macaulay.\end{corollary}

        \begin{proof} From Theorem~\ref{main theorem} and
        Proposition~\ref{F1} it follows that $I$ satisfies sliding
        depth and the condition $\F _1$, and so $R[It]$ and ${\rm
        gr}_I(R)$ are Cohen-Macaulay by Theorem~4.3.8 and
        Corollary~3.3.21 of [V]. 

        Proposition~\ref{F1} and Theorem~\ref{main theorem} also imply
        that for all $t$, the modules $I^t/I^{t+1}$ are torsion-free
        over $R/I$ (using approximation complexes), and since $I$ is
        generically complete intersection and $R/I$ is reduced, it
        follows from Corollary~5.3 of [SVV] that $R[It]$ is normal.
        \end{proof}

%%%%%%%%%%%%%%%%%%%%%%%%%%%%%%%%%%%%%%%%%%%%%%%%%%%%%%%%%%%%%%%%%%%%%%%%%%%%%

{\bf Acknowledgments.} I would like to thank Wolmer Vasconcelos for
reading an earlier version of this manuscript and for many helpful
comments. I am grateful to Dan Ullman and Peter Selinger for
illuminating conversations on trees, and especially to Will Traves,
for introducing me to edge ideals and for many discussions on facet
ideals during the early stages of this work.

%%%%%%%%%%%%%%%%%%%%%%%%%%%%%%%%%%%%%%%%%%%%%%%%%%%%%%%%%%%%%%%%%%%%%%%%%%%%%

\end{document}